\RequirePackage{fix-cm}

\documentclass[smallextended]{svjour3}       

\smartqed  

\usepackage{graphicx}
\usepackage{epstopdf}

\usepackage{amsmath}
\usepackage{mathrsfs}

\usepackage{algorithm}
\usepackage{algorithmic}

\usepackage{booktabs}
\usepackage{enumitem}

\begin{document}

\title{Multiscale Stochastic Simulation of the US Pacific Northwest Using Distributed Computing and Databases with Integrated Inflow and Variable Renewable Energy}
\subtitle{The new Genesys model}


\author{Joaquim Dias Garcia         \and
        Guilherme Machado         \and
        Andr\'e Dias    \and
        Gerson Couto \and
        John Ollis \and
        John Fazio \and
        Daniel Hua
}


\institute{Joaquim Dias Garcia, Guilherme Machado, Andr\'e Dias and Gerson Couto are with PSR at:\at
              Praia de Botafogo, 370 - Botafogo, Rio de Janeiro - RJ, 22250-040, Brazil \\
              Tel.: +55-21-3906-2100\\
              \email{\{joaquim, guilherme, dias, gerson\}@psr-inc.com}           
                                          \\
              John Ollis, John Fazio and Daniel Hua are with Northwest Power and Conservation Council, Portland, USA\\ 
              \email{\{jollis, jfazio, dhua\}@nwcouncil.org}
}

\date{August 5, 2020}

\maketitle
%


\begin{abstract}
Modelling challenges of the United States Pacific Northwest system have grown in the last decade. Besides classical modelling difficulties such as a complex hydro cascade with many operational constraints, we have seen higher penetration of variable renewable energy inside and outside the system leading to internal issues and completely different power exchanges with the West Coast System. The analysis of adequacy and reliability of this system motivated the design and implementation of a five-step simulator including four planning phases and an operation step. The five levels were modeled as mathematical programs that are linked from top to bottom by fixed decisions and improvements in forecasts, on the other hand they are also linked from bottom to top by system updated states. The solution of the millions of resulting mathematical programs was made possible by applying state of the art optimization techniques and high performance data bases in a massively parallel environment.
\end{abstract}



\section{Introduction}

Analytical modeling of the 63.5-GW US Pacific Northwest (USPN) has historically been challenging because of the complex Columbia river operation rules for flood control, Canadian upstream storage, salmon management and many others. In the past years, this complexity has been compounded by the integration of gas-fired plants and the massive construction of variable renewable energy (VRE) such as wind and solar (15\% of total capacity), which led to higher variability and transmission issues. Another disruptive change has been the very high penetration of solar power in California and other states in the West Coast system (WECC), which affected the pattern of market supply and economics of power exchanges between the USPN and those states.

In response to those new modeling challenges, the Northwest Power and Conservation Council (NWPCC) decided to develop a new version of their stochastic operation model, Genesys. The main objective of the new Genesys is to represent in detail the hierarchical sequence of operational decisions, composed of five planning/scheduling phases: (i) year-ahead planning; (ii) week-ahead planning/scheduling; (iii) day ahead scheduling; (iv) hour ahead rescheduling; and (v) real time operational decisions ("true-up"). In each phase operational decisions are taken under different degrees of forecasting accuracy (load, wind, solar and inflows) and updated information on system conditions (e.g. equipment outages).

NWPCC specified that steps (i)-(v) should be carried out/updated chronologically for each of the 8760 hours of the year, in order to capture the sequential nature of operational decisions. Additionally, the entire 8760-hour process should be replicated thousands of times for different inflow, renewable, load and outage scenarios. The new Genesys outputs include supply adequacy measures and probability distributions of operational indices such as hydro and thermal generation, fuel consumption, power exchanges with WECC and others.

It will be seen that the above simulations require the execution of up to a hundred million large-scale Mixed Integer Programs (MIP) and Linear Programs (LP) optimization problems.
Moreover, the total execution time for a complete study should be limited to allow sensitivity studies, for instance, eight hours, as in the original specification.
Finally, because the model was intended to be used by many stakeholders, ranging from government agencies (e.g. US Army Corps of Engineers), utilities (e.g. BPA), private investors, environmental groups and others, NWPCC indicated that open source tools should be used as much as possible.

This paper describes the methodologies, analytical tools and computational resources integrated into the new Genesys system: (i) multistage stochastic, MIP and affine optimization for the planning/scheduling models; (ii) Bayesian networks for the production of integrated multiscale inflow and renewable scenarios; (iii) representation of forecasts as weights on those scenarios, which change with the operational horizons; (iv) cloud-based distributed computation (more than 300 processes, each with 36 cores); (v) cloud-based distributed data management resources (Cassandra, Amazon S3, Spark and Amazon Athena) to handle the very large files (+10 Tb) of hourly renewable scenarios and result outputs; and (vii) open source optimization languages MOSEL and JuMP (a Julia library).

Several new simulation tools for power system planning have been proposed in the past few years. In addition to proprietary software, we highlight open source tools such as \cite{chassin2008gridlab} and \cite{anderson2014gridspice}, which were specially designed to smart grids. Simulators for integrated transmission and distributions software are detailed in  \cite{huang2016integrated}; efficient tools for unit commitment studies are presented in \cite{fu2007fast}. Systems with large penetration of water resources often use a hierarchical approach integrating models with longer- and shorter-dynamics, see e.g.
\cite{belsnes2016applying}. Hierarchical systems were also studied in: \cite{ilic2012hierarchical} with focus on higher frequency dynamics from unit commitment to Automatic Generation Control (AGC); and in  \cite{Operations2019Atakan} with focus on high renewable energy penetrations.

Section \ref{details} provides an overview of the multiscale simulation approach and the associated computational challenges. The methodology is detailed in section \ref{method}, and the computational architecture is described in section \ref{arch}. In section \ref{study} we present a case study to illustrate the model performance and capabilities. Finally, section \ref{conclusion} presents some conclusions and future work.

\section{Multiscale power system simulation}\label{details} 

Because of the combination of long-term horizons, multiple uncertainties (loads, inflows, renewable production, fuel prices, equipment availability and others) and the need for hourly (or smaller) resolution, it would be both cumbersome and computationally infeasible to model power systems operation as a single mathematical optimization problem. Therefore, a multiscale hierarchical approach is widely used in actual systems. First, "slow frequency" dynamics, e.g. hydro opportunity cost, are made, followed by also "slow" unit commitment decisions (plants with long up-times and slow ramps). Finally, as we get closer to the actual operation, fast thermal plants, storage and other dispatchable resources are (re)scheduled in order to handle VRE and demand fluctuations.

A hierarchical scheme that mimics many steps of real-life operation of a power system is key to evaluate reliability and adequacy, because the physics and existing rules for the system operation can be tested against many hypotheses ranging from past historical data to forecasts of climate change scenarios. Detailed realistic simulations can provide valuable information about the behavior and weaknesses of the system in question, enabling the decision maker to plan in advance, determine good policies and avoid stress situations.

Although the new Genesys model was tailored to the USPN system requirements, the methodology, analytical tools and computational architecture is general enough to be applicable to other systems. The reason is that, as mentioned, most power systems use hierarchical scheduling and operation steps, in which optimization problems are solved considering past decisions, current system state and forecasts. The simulation proceeds from one layer to the other in a sequential way, trying to mimic the real-life operation.

\subsection{Operation layers: decisions, forecasts and details}

\begin{figure}[htb]
\centerline{\includegraphics[]{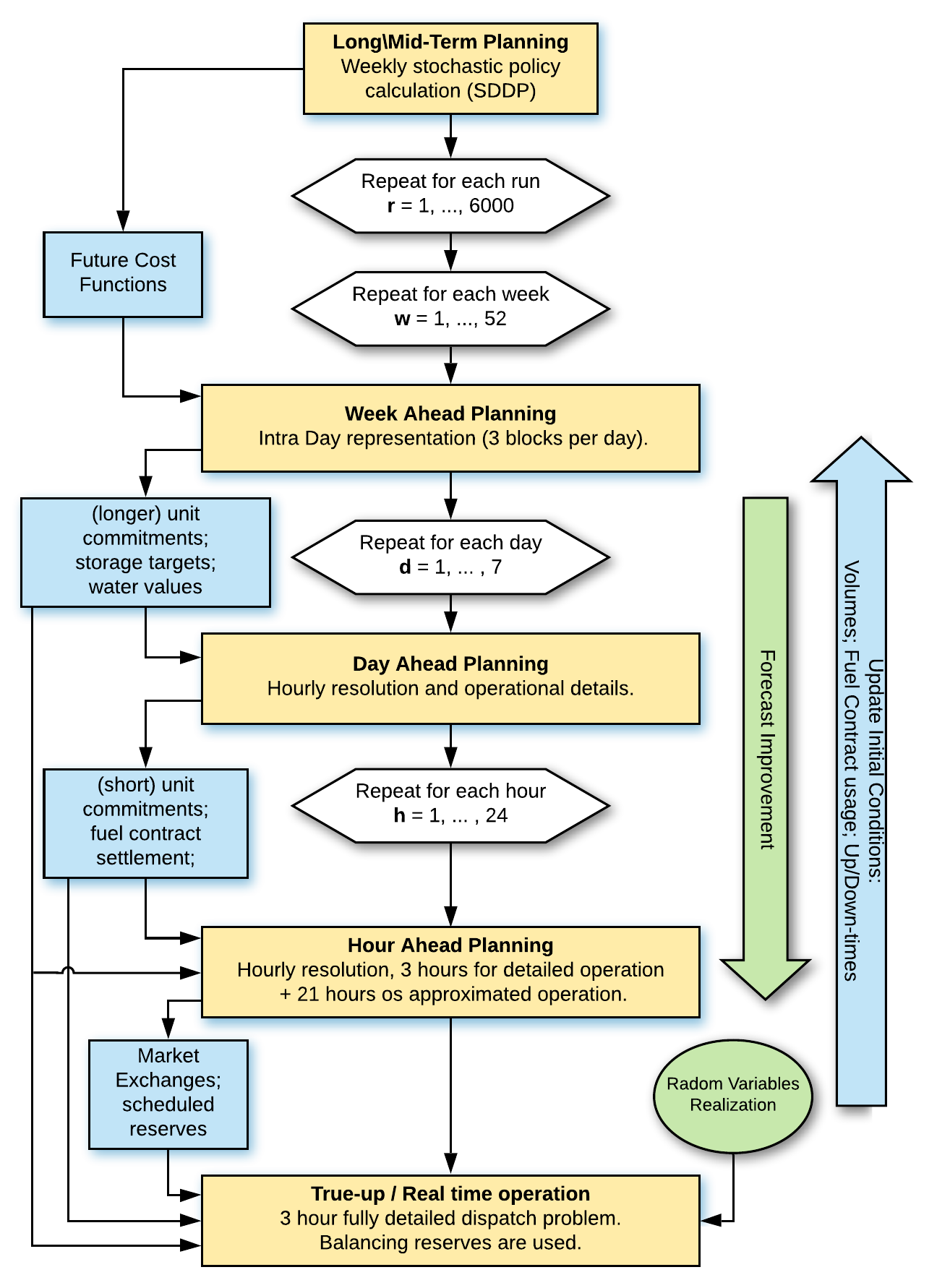}}
\caption{Multi level model with information flow}
\label{fig:multilevel}
\end{figure}

The five layers in the simulator are: Mid(or long)-term planning, a optimization model with multiple weeks as optimization horizon; week ahead planning; day ahead planning; hour ahead planning; and real time operation, referred here as true-up phase. Fig. \ref{fig:multilevel} presents all the five levels, depicted in the central yellow boxes with bold titles, other important elements of the simulation are highlighted. White hexagons represent the repetitions of the many planning and operation phases. Blue boxes on the left side of the figure describe information that is being passed to from one step to the other. The green down pointing arrow highlights that at each level the forecast quality is improves until the actual realization of random variables in the true-up phase. Finally, the up pointing blue arrow emphasizes that the system state is updated in the above levels after the new system state is reached in a lower level.
In the next paragraphs we detail the main elements of the simulator present in Fig. \ref{fig:multilevel}.

\subsubsection{Mid-term planning}

The first layer is responsible for strategic long \slash mid-term decisions. The slow dynamics are the focus of the optimization, fast dynamics are also considered, but with less details. This step usually encompasses multiple months, weeks or even a few years. Currently, this layer is responsible for optimizing the reservoir usage policy or defining water opportunity costs, because inflows processes and large reservoirs inherently require this kind of mid-term time coupling. Because of the large time scale, decisions are made under uncertainty by formulating the problem as a Multi-stage stochastic program (MSP) \cite{pereira1989optimal}. Many hydro constraints are considered, for instance flood control water levels, bounds on spillage, energy production and other reservoir level constraints. VRE and demands are represented as scenarios, that is a realization of random variables.

\subsubsection{Week ahead planning} With the optimal mid-term policy determined, in our case represented by water value curves (future cost functions), the simulation proceeds to weekly decisions. In this phase, almost nothing is known about VRE and demands and therefore very naive forecasts are considered. Decisions of this time step include slow commitment plants (with up-times longer than days), reservoir level targets or detailed water opportunity costs.

\subsubsection{Day ahead planning} At this point, forecasts are improved and this new information is combined with previous decisions to optimize a 24 hours detailed planning problem. Decisions taken at this point include reserves scheduling, commitment of plants with intermediate up-times from 3 to 10 hours and gas contract nominations. These decisions are fixed and the simulation proceeds to the next step.

\subsubsection{Hour ahead planning} In the last of the planning steps, many of the system resources are already allocated, but new and more accurate forecasts are made available. This allows the system operator to review some decisions, determine commitments of fast plants, decide market exchanges at updated market prices. At the end of this step the system is ready for real time operation with limited available resources.

\subsubsection{True-up phase} Finally, the real VRE and demand is observed and the systems must adapt to withstand these changes. Since all the commitments and the area exchanges are fixed, each area must operate the previously allocated reserves to balance load and generation. At this final step it is also possible to consider events like generator and line failures, which are accommodated by dispatching contingency reserves.

\subsubsection{Sequential decisions}

In actual planning and operation of power systems the results of the last planning phase and the current state of the system should be fed into the next steps. Therefore, communication between layers also go from the bottom to the top. The last true-up reveals a power system state (reservoir levels, fuel contract availability, and thermal up and down times) that must be fed into the next hour ahead planning step. The same happens when solving the next day and week ahead planning problems that require the current state of the system.

In other words, steps (i)-(v) are carried out until the true-up of 1 AM of January 1st; next, all forecasts and operating conditions for 2 AM and following hours are updated based on the observed values of hour 1; the process is then repeated until the true up for 2 AM; and so on until hour 24 of December 31st.

\subsection{Computational Challenges}

Simulating all these layers is a complex problem due to the modelling details. The very first operation layer is formulated as an MSP, which is an extremely hard class of problems, even for convex sub-problems. Therefore, special algorithms must be applied.

The three layers ranging from week ahead to hour ahead require the solution of complex problems firstly due to the non-linearity in hydro production factors, hydro constraints, integer nature of the commitment decisions and the time coupling of up and down times, water mass balance and gas nomination contracts. Although these layers are not optimized under uncertainty, the number or problems that are solved is huge because we are simulating the system for many scenarios.

Finally, the true-up phase must contain another hourly simulation with the realization of the uncertainties and the final dispatch of the system's resources. At this point, if the planning occurred as expected, the system must be able to supply the demand considering the commitment decision of previous phases, fueling contracts and assigned reserves.

\section{Methodology}\label{method}

The simulation scheme detailed in the previous section is a multi-level hierarchical model. Each of the layers is formulated using different optimization frameworks and solved by specialized methods. Not only are the problems being solved intrinsically complex, but also the number of scenarios can be very large. Therefore, the models formulated in each of the simulation layers must be solved by specialized methods. We begin by describing large-scale data required for the simulation, then we detail the solution method applied in each layer and we finish with a brief description of the forecasting simulations that relies solely on scenario data.

\subsection{Input data}

Even before the first level the simulator requires data to be prepared. Some of the data are related to resources (physical data, policy rules, economic and biological constraints). Data is usually small and sparse, and in some of the most complex systems can reach a few hundred megabytes. The remainder is time series related data: variable renewable energy scenarios like wind and solar production; hydro inflows; inelastic demand and demand response scenarios. This data can be massive because we need values for all the hours of the year, and we have as many years as scenarios.

Scenarios usually come from historical data: either the historical data itself is used to construct scenarios or time series are simulated by chronological Monte Carlo methods. Modelling large basins, multiple VRE sites and considering wind and inflow correlation is a complex task as presented in \cite{souto2014high}. We propose to use a Bayesian network approach to estimate spatial correlations combined with a periodic auto regressive model to capture temporal correlation while generating time series data \cite{TSL}.
The correlations between such random variables can either be estimated directly or considering latent variables such as temperature. The large dimension of the resulting data will require specialized infra-structure.

\subsection{Mid-term planning}

In the very first layer, bulk data is used to formulate a Multi-Stage Stochastic Optimization problem  \cite{shapiro2009lectures}. This problem is responsible for determining the optimal allocation of resources with yearly dynamics, for example water in reservoirs. Each stage of this problem corresponds to a single week in a total of $T$ weeks, one year typically, and the number of scenarios, $S$ is given as input. Weekly stages contain 3 chronological load blocks per day, instead of 168 hours, to keep the problem tractable.

Instead of regular scenario trees,
frequently used in stochastic programming literature we use samples of a scenario tree or paths in the scenario tree because they fit perfectly with the solution algorithm employed: Stochastic Dual Dynamic Programming \cite{pereira1991multi} SDDP. SDDP is a state-of-the-art method for MSO because it scales well with the number of scenarios and state variables, such as reservoir levels. 

The solution of the MSO problem in the first layer generates a weekly policy encoded in the SDDP's future cost function (FCF), which indicates the average value of resources of yearly dynamic, such as water, in the current week. This is a common part of the Brazilian official planning process \cite{maceiral2018twenty} and a widely used techniques in countries with relevant water reservoirs such as Norway \cite{rotting1992stochastic}, Canada \cite{meier2012real} and New Zealand \cite{halliburton2004optimal}. The FCF is then passed to the following layers to complete the simulation.

\subsection{Week, Day and Hour ahead planning}

These three planning steps have a lot in common. All of them require the solution of Mixed Integer Programs (MIP) with many details of the system such as commitment data (up-time, downtime, ramps) and non-linear hydro production functions.

On the one hand, week ahead planning is solved with 21 steps just like the Mid-term planning and it is responsible for converting the FCF into hourly water values. These values are then passed to the following layers together with reservoir daily level targets. Long commitments, for plants with long up and down times, are also decided at this point. All this optimization is processed with limited knowledge of the random variable that will be observed during the last time step, only a, typically poor, forecast is available.

On the other hand, the day and hour ahead planning are both carried out with daily horizons and (24) hourly steps. In the day ahead planning the forecasts are improved and the main results are thermal plant commitment decisions, gas contract amounts and allocation of reserves.

Finally, in the hour ahead planning, most of the commitments have already been fixed and now we have a target for gas that might be used as planned or not, because the forecast is better for this time step. The main decisions here are: some very fast plants commitments, demand response and market transactions. The last-minute changes here are crucial because in the real time operation each area will have to be operated individually by using the previously allocated reserves.

Note that the hour ahead problem is a 24-hour problem solved in a rolling horizon framework. This framework is necessary because approximations are improved if there is some extra information about the future since there is some time coupling in reservoir levels and thermal commitments.

The key aspect in this section is to highlight some crucial techniques employed in these MIPs. The most important is the application of strong unit commitment formulations \cite{ostrowski2011tight} \cite{guan2018polynomial} and the second is to reduce the number of integer variables, not all hours of the hour ahead planning problem require representation as integer since they are only useful to improve approximations of the unknown future. 

\subsection{True-up phase}

The true up is the last phase and it is responsible for the final evaluation of all the planning process against actual observed realizations of the random variables. In this phase, it is solved a 3 hour problem, the first hour corresponds to actual realizations of the random variables while the following are forecasts. This phase is computed for each hour of the simulation. Since it uses an hourly representation, we consider a stochastic formulation, thus we can better represent the uncertainty that could happen within the phase. Inflow scenarios are fixed at this point, but we have weighted (not uniform) distributed demand and VRE scenarios.

In order to be able to solve stochastic programs in reasonable time we rely on a method called affine decision rules\cite{Bodur2018a}. This method consists in formulating the problem so that scenario dependent decision variables are converted to functions that depend linearly on the random variable, and the parameters of this functions are optimized instead of the individual scenario decisions. This drastically reduces the number of decision variables in the optimization problem and makes it computationally tractable \cite{Kuhn2011}\cite{Shapiro2006}, at the expense of giving a sub-optimal \cite{Shapiro2006} or infeasible \cite{Chen2008} result, although \cite{Lorca2016a} showed that a simplified affine policy is quite powerful and produce close-to-optimal results for a robust multistage unit-commitment problem. In power systems, this method also has an interesting interpretation: the decision rules can be interpreted as automatic generation controls (AGC), properly optimized to operate at a given set-point but also prepared to adjust the operation in face of deviations from the set-point \cite{Jabr2013}. The linear decision rules can be interpreted as the participation factors that determine the change in the power generation due to a given deviation in the scheduled production.

The final solution is the one tied to the scenario corresponding to the actual realization of the random variable.

\subsection{Forecasting}\label{sec:forecasting}

A key feature of this hierarchical system is the ability to start with a given set of time series scenarios and be able to emulate the forecasts used in planning steps. In this section we detail how this works.

Given a set of $S$, not weighted, realizations for some random variable we choose one of them to be the actual realization $s$. After that we must define weights for all these scenarios, considering the that forecasts have no or small bias, we set weights considering that the largest weight is given to the $s$ while the other weights are increasingly smaller for scenarios more distant from $s$. One way of doing so is to follow a weighted distribution such as a Gaussian: we set the scenario $s$ as the mean and set the variance (of the Gaussian) in order to have a random variable with a given variance (of the weighted sample), chosen by the decision maker.

For each layer it is possible to define a different target variance, or another metric such as Coefficient of Variation (CV) and Root Mean Square Error (RMSE). By computing the weighted average of the scenarios, we have a forecast of that random variable corresponding to that realization. As the simulation layer is closer to the final evaluation the variance is set to a lower value, and the forecast is improved. Doing that for all the hours of the study we can obtain one time series for each forecast level for a given scenario $s$, for each layer.
We repeat the same process for all the scenarios $S$ and obtain $S$ time series for each layer of the study.

These important features allow the simulation to have higher fidelity and avoid over fitting in the initial planning phases. Moreover, this allows performing sensitivity analysis of good forecasts and its effects in the overall system operation. Finally, the sensitivity analysis can be used to assess the value of forecast improvements.

\section{Computational architecture}\label{arch}

Simulating large-scale power systems for each of the above-mentioned decision phases is a difficult task on its own. Doing that for all the phases in sequence with many scenarios is an overwhelming time and resource consuming task. In order to perform the complete system simulation in reasonable time, say 8 hours, we need not only to use state of the art optimization method such as SDDP and first class optimization libraries for LP and MIP solving, but also code must be written in high-performance frameworks with great usage of parallel computing, and high performance IO (input output) for which data bases are required. In usual applications, optimization solvers are the bottleneck of the solution processes, however, writing all detailed solutions for all the hours in each layer, for all scenarios, accounts for more than 10 TB of data, which cannot be stored in regular systems.

\begin{figure}[htb]
\centerline{\includegraphics[width=1\linewidth]{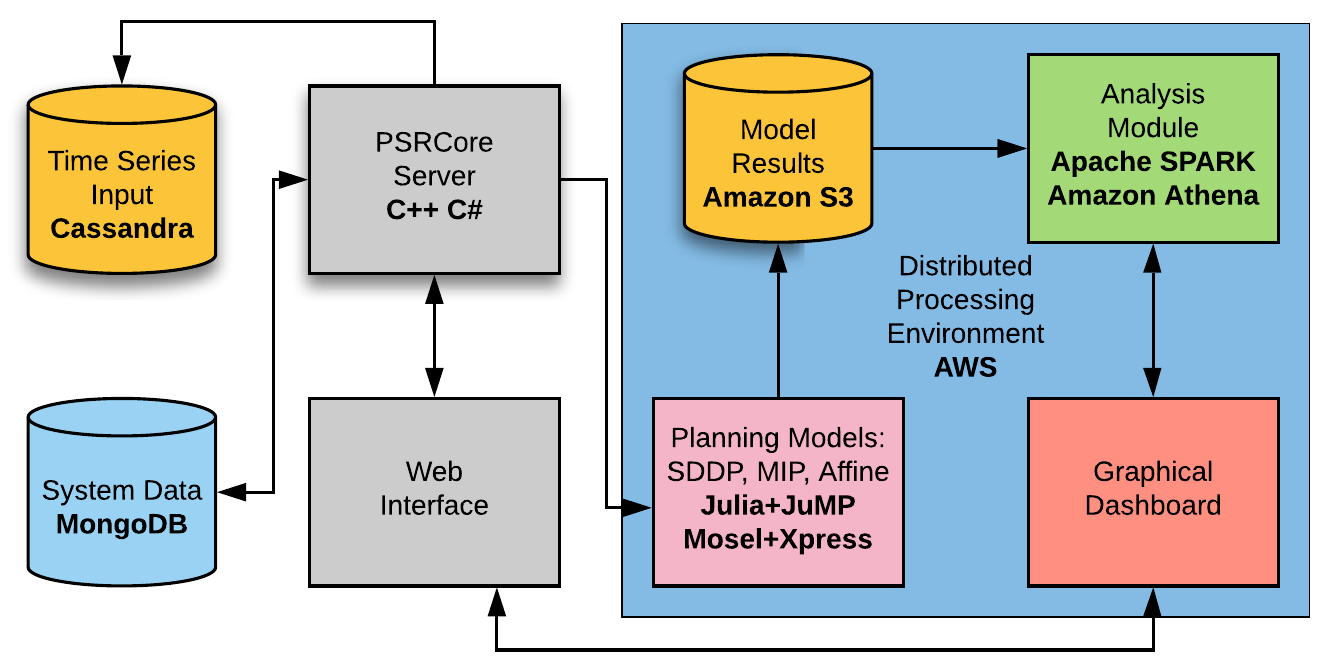}}
\caption{System architecture}
\label{fig:arch}
\end{figure}

The complete architecture is described in Fig. \ref{fig:arch}. In the leftmost column we find the data bases used to store input data, the following column highlights the central server responsible for passing information to the web interface (right below the server box) and to the multi-level model (pink box). The large box on the right side is a massively parallel environment where the models are executed with access to high performance data bases also where results can be queries with analysis and graphical dashboard tool. The main block, distributed computing, will be detailed in the next paragraphs.

\subsection{Parallel computing}

Massive parallel computing can be used in all steps of the solution approach. In the first layer, we have the SDDP method that is very amenable to parallel computation \cite{da2003parallel}. Although it is an iterative algorithm, all scenarios are dealt with in each iteration and that is done in parallel.

The three following planning phases depend of the solution of MIPs, for each scenario. In each layer we solve problems for all scenarios at the same time because there is no communication between them and therefore a few computing cores can be used in each scenario because the solver, using a few cores, speeds up solution times.

The final true-up layer can also be solved the same as the previous three because the true-up of one scenario will not affect the true up of another one. In fact, for each scenario we solve a true-up problem considering multiple scenarios, but all of them have different weights, because each one represents the real-life operation for one possible realization.
This massive parallel computing can be performed in hundreds or thousands of cores. Instead of relying on a on-premise infra-structure for this computation we can use cloud computing such as Amazon Web Services (AWS)  \cite{jackson2010performance}, with performance tuned computing cores. One extra advantage of this design is that it is possible to adapt the number of servers used in each phase and the number of servers used in each simulation. If some study requires fewer scenarios and time steps, a smaller and consequently cheaper infra-structure can be applied.

\subsection{High performance data bases}

High performance data bases are important for at least three fundamental steps in the simulation: reading data, writing data and querying results.
In a massive parallel application, reading input data and writing solutions become relevant because usual cloud computers have IO limits, therefore we designed a system to rely on special databases.

For the case of time series input data, we have chosen Cassandra \cite{chebotko2015big}, a high-performance data base. Cassandra allows the system to read data in parallel without the access problems of standard data files.

While designing the systems to  write solutions two features must be considered: first, the processes of actually writing the outputs must not be a bottleneck of the system, performance is fundamental; second, the saved solutions must be readily and easily accessible for the decision maker to perform analysis and extract information. The tool selected was Amazon S3 \cite{palankar2008amazon} because it is inherently tied to the amazon cloud, it has good performance and is a relatively cheap data storage. Moreover, S3 integrates seamlessly with the Amazon Athena \cite{documentation2018amazon} service that can perform queries in extremely large data bases in a very reduced time, solving the problem of making the results of the simulation easily available for the decision maker.

\subsection{Optimization software}

Finally, all the parallel implementation must be orchestrated by software and, in this simulator, we used Julia \cite{bezanson2017julia}, a high-performance language for scientific computing. Besides being extremely productive for writing parallel software and high performance code, Julia also has very good libraries for operations research and optimization; we highlight JuMP \cite{DunningHuchetteLubin2017}, a mathematical programming language used for coding our SDDP algorithm. References of interesting applications of Julia and JuMP in power systems can be found in \cite{8442948}; a open source SDDP implementation in Julia can be found in  \cite{dowson_sddp.jl}.
The other optimization problems were written in MOSEL  \cite{colombani2002mosel}, which has recently been made open to the general public; it is very efficient and makes the process of problem formulation and tuning the MIP solver very productive.

\section{Case Study}\label{study}

In this study, the USPN system is represented in detail. We have individual representation for the main American and Canadian hydro plants in the Columbia River.
The list of plants includes 76 hydro plants: 38 reservoirs and 38 run-of-river plants. The main thermal energy resources of the detailed represented system include 133 thermal plants, 51 of which have fuel contract constraints and 90 have commitment constraints. The remainder of the generation units are represented as non-dispatchable units: 127 wind power plants, 30 solar power plants, 54 small hydro plants and 128 independent resources.

Although the simulator can represent electric networks in detail, in this case study we consider a simplified network. The network includes some of the most important buses in the system and some other buses to which the system is connected: this adds up to 34 buses with individualized hourly loads with inelastic and elastic segments. Completing the simple network, we have 71 circuits connecting these buses.

In many systems the market trades can be as important as the system's own resources, for USPN we represent 18 buses with access to energy market represented by price and quantity curves. The average number of segments in those price and quantity curves is 21.
Balancing areas are important regulatory division in the region. The simulator represents 21 of them including all their resources (including shared resources) and individual reserve constraints that are used to perform redispatch in the system due to divergences from forecasts to actual realizations of the VRE and demand random variables.

We were able to perform the complete simulation of one year in 7.5 hours. The SDDP algorithm in the first step takes around 2000 seconds to complete the simulation. After that the systems have to process all the scenarios independently, each one of them requires the solution of almost 18 thousand MIPs, two for every hour of the year (hour ahead and true-up), one for each day and one for each week. This second step requires 25 thousand seconds (almost 7 hours) to complete for each scenario. By deploying one core per scenario we can finish the computations in the aforementioned time. The EC2 machines used have 36 cores each, if we allocate 2 cores per scenario, we can complete a 6 thousand scenario simulation with 334 servers.

\begin{figure}[htb]
\center
\includegraphics[trim={0 0 0 2cm},clip,width=1\linewidth]{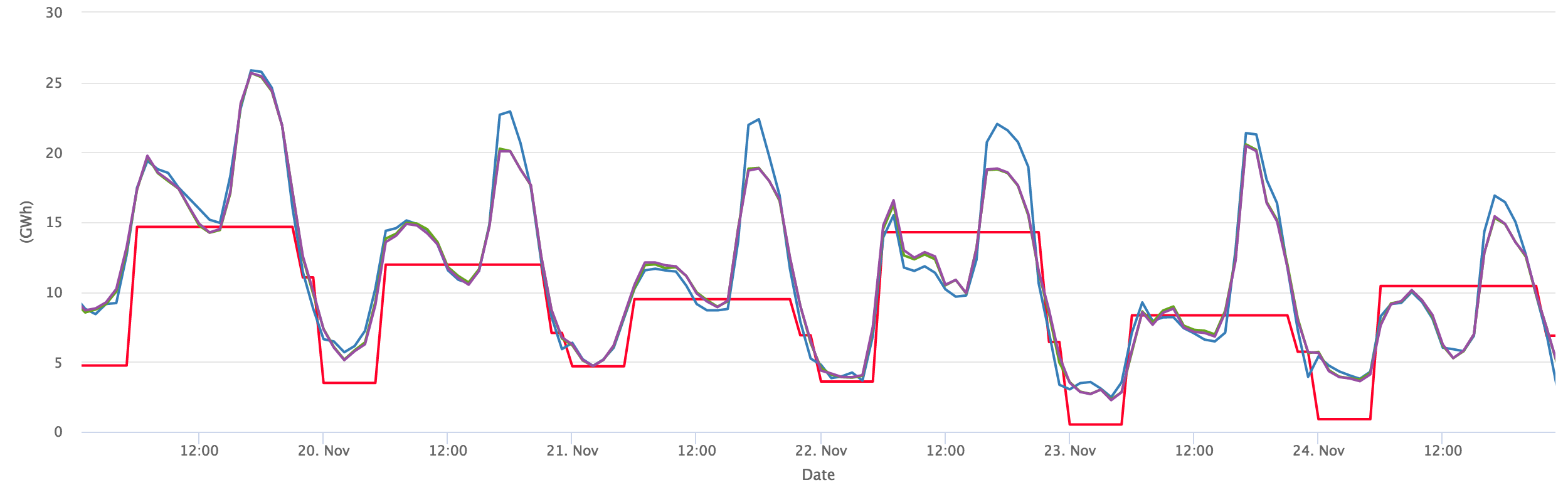}
\caption{Hydro generation, six days. Week ahead (red),
day ahead (blue), hour ahead (green) and true-up (purple).}
\label{fig:hyd_week}
\end{figure}


\begin{figure}[htb]
\center
\includegraphics[trim={0 0 0 5cm},clip,width=1\linewidth]{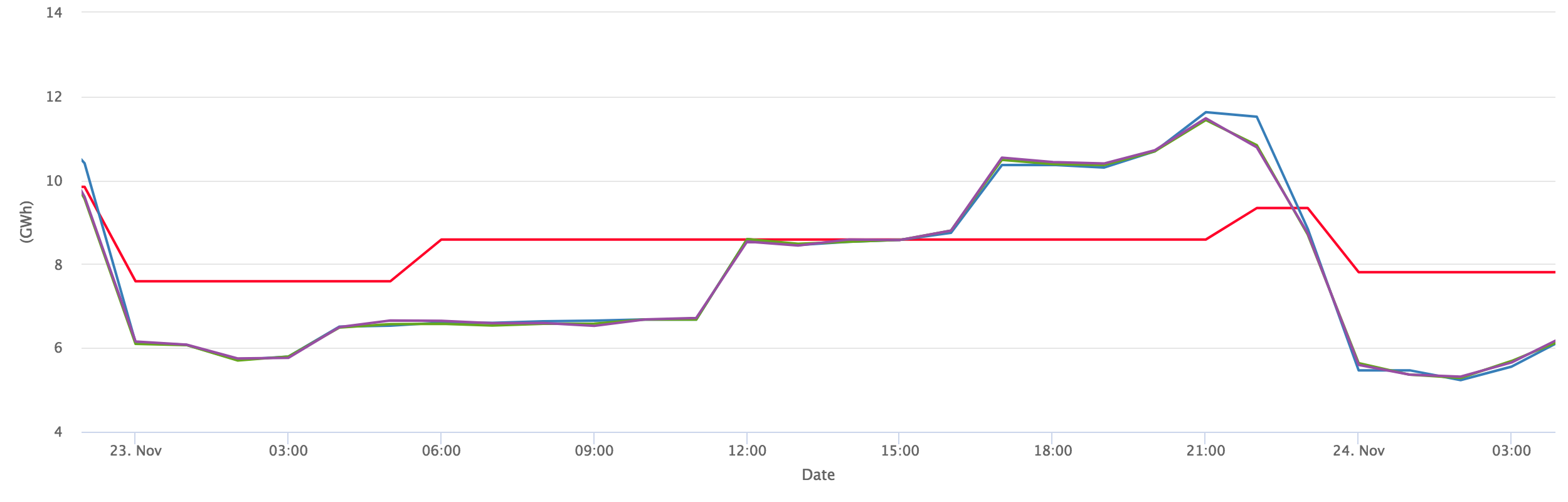}
\caption{Thermal generation, single day, good forecast. Week ahead (red),
day ahead (blue), hour ahead (green) and true-up (purple).}
\label{fig:ter_day_good}
\end{figure}

\begin{figure}[htb]
\center
\includegraphics[trim={0 0 0 3cm},clip,width=1\linewidth]{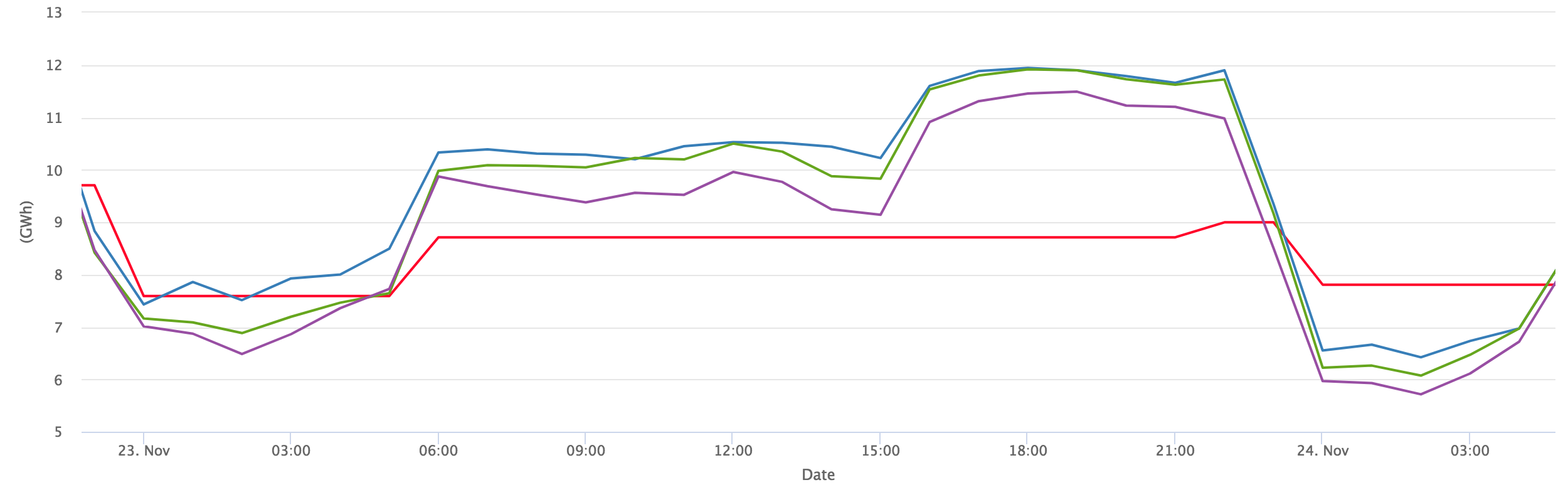}
\caption{Thermal generation, single day, poor forecast. Week ahead (red),
day ahead (blue), hour ahead (green) and true-up (purple).}
\label{fig:ter_day_poor}
\end{figure}

As an example of the simulation performed by the system we present results for one fixed scenario: Fig. \ref{fig:hyd_week} shows six days of hydro operation, Fig. \ref{fig:ter_day_good} displays one day of thermal generation and Fig. \ref{fig:ter_day_poor} presents the same day, but the operation is done with a much worse forecast. In the three figures we have four phases, the first one is only represented by FCFs and, therefore, not presented.

We can see in Fig. \ref{fig:hyd_week} that the raw decisions of the week ahead plan is coherent with the other decisions, which are initially off in the day ahead but approaches the final true-up operation.
Contrasting Figs. \ref{fig:ter_day_good} and \ref{fig:ter_day_poor} we can see the importance of analysing forecasting quality: the mismatch between plan and reality increases together with the thermal generation (and  cost) because the planned decisions will more distance from the optimal ones in cases of poorer forecasts.

\section{Conclusions}\label{conclusion}

In this work we proposed a simulation model, motivated by the need of evaluating system adequacy and reliability in a world that gets more complex every day. The simulator is very well suited to the USPN system, which has complex hydro system, growing penetration of VRE and connections to the entire WECC and Canada. A five-level model was proposed with planning phases: year, week, day and hour ahead; and an operation evaluation phase.
This model is solved by carefully modelling all the layer within the mathematical programming paradigm. Solution techniques were carefully chosen and implemented, we highlight SDDP, Affine Rules and MIPs. In order to represent operation of systems with relevant VRE a progressive forecasting emulator was proposed.

The simulator turned to be an extremely computation intensive program that required a tailor-made architecture. We combined modern computing languages and data bases in a massively parallel environment. This allowed the multi-step simulation of a fairly detailed representation of the USPN in reasonable time.

Finally, we highlight next step in the research. The representation of more complex networks is a important feature that is already under testing. More computationally intensive systems, such as the USPN including more detailed representation of resources in the surrounding markets, will also be tested in future work.



\bibliographystyle{IEEEtran}
\bibliography{ref.bib}
%



\end{document}